\documentclass{amsart}
\usepackage{amssymb,latexsym,verbatim,hyperref}
\usepackage[all]{xy}

\newcommand\DMO[2]{\DeclareMathOperator{#1}{#2}}
\newcommand\comp{\circ}

\newcommand\Rn{\mathbb{R}}

\newcommand\qbound{\preccurlyeq}
\newcommand\qequiv{\approx}

\newcommand\homot{\simeq}
\newcommand\bdy{\partial}
\renewcommand\bar[1]{\overline{#1}}
\newtheorem{theorem}{Theorem}
\newtheorem{corollary}{Corollary}
\newtheorem{lemma}{Lemma}
\theoremstyle{definition}
\newtheorem{definition}{Definition}

\DMO{\Obj}{Obj} \DMO{\Mor}{Mor} \DMO{\im}{im} \DMO{\coker}{coker}
\DMO{\Hom}{Hom} \DMO{\Spec}{Spec} \DMO{\Frac}{Frac} \DMO{\tr}{tr}
\DMO{\Sym}{Sym} \DMO{\supp}{supp} \DMO{\Proj}{Proj} \DMO{\vol}{vol}
\DMO{\Riem}{Riem} \DMO{\Diff}{Diff} \DMO{\FV}{FV} \DMO{\Vol}{Vol}
\DMO{\dist}{dist} \DMO{\Cont}{Cont}

\newcommand\Volcur{\mathbf{M}}
\newcommand\Volmet{\Vol_\mathrm{met}}

\newcommand\Voltr{\Vol_\mathrm{tr}}

\newcommand\FVcur{\FV_\mathrm{\!cur}}
\newcommand\FVmet{\FV_\mathrm{\!met}}
\newcommand\FVch{\FV_\mathrm{\!ch}}
\newcommand\FVcell{\FV_\mathrm{\!cell}}

\newcommand\Phicur{\Phi_\mathrm{cur}}
\newcommand\Phimet{\Phi_\mathrm{met}}
\newcommand\Phich{\Phi_\mathrm{ch}}
\newcommand\Phicell{\Phi_\mathrm{cell}}

\renewcommand\tilde{\widetilde}

\begin{document}

\title{Generalized Dehn functions II}
\author{C. L. Groft}
\date{January 15, 2009}

\begin{abstract}
  For $G$ a group of type $F_q$, we establish the existence,
  finiteness, and uniqueness up to scaling of various $q$-dimensional
  isoperimetric profiles. We also show that these profiles all coincide
  for $q\ge 4$, and that significant overlap exists for $q=3$. When $G$
  has decidable word problem, this has mild consequences for the growth
  rates of these profiles. We also establish a metric analogue for highly
  connected Riemannian manifolds.
\end{abstract}
\maketitle

\section*{Introduction}\label{S:intro}

In part I (see \cite{clG08}) we examined the isoperimetric profiles
$\Phi^{X,M}$ where $M$ is an orientable manifold of dimension $q\ge 2$ with
nonempty boundary and $X$ is either a CW complex or a local
Lipschitz neighborhood retract (LLNR). Briefly, if $f\colon\bdy M\to X$
is a map with volume $v$ that extends to $M$, then $f$ extends
to $g\colon M\to X$ with volume at most $\Phi^{X,M}(v)$, and
$\Phi^{X,M}(v)$ is the smallest nonnegative extended real number
with this property. In order for a map $g\colon M\to X$ to have a volume,
it must be a map $(M,\bdy M)\to(X^{(q)},X^{(q-1)})$ (which in \cite{clG08} is called
\emph{quasi-cellular}); there are several possible
definitions in this case, as discussed in \cite{clG08}. When $M=D^q$, these profiles are identical
to the higher-dimensional Dehn functions defined in \cite{jmA99}.
Similar profiles $\Phi^{X,q}$, for which we replace
maps $\bdy M\to X$ and $M\to X$ with $(q-1)$- and $q$-dimensional chains
respectively, were also defined, both in \cite{clG08} and earlier in
\cite{dbaE92} and \cite{mG93}.

Unlike the more familiar Dehn function of a complex (think $M=D^2$),
the functions $\Phi^{\tilde X,M}$ and $\Phi^{\tilde X,q}$ (for $X$ a finite complex)
are not dependent on $\pi_1(X)$ alone; the homotopy groups $\pi_2(X)$ through
$\pi_{q-1}(X)$ are also relevant. Moreover, if the higher homotopy groups are
nonzero, it is possible for $\Phi^{\tilde X,M}$ and $\Phi^{\tilde X,q}$ to be different
functions, and $\Phi^{\tilde X,M}$ may depend nontrivially on $M$ (this is not yet
clear). The best theorem we have at present is the following: Given a continuous
function $f\colon X\to Y$ where $f_*\colon \pi_t(X)\to\pi_t(Y)$ is an isomorphism
for $1\le t <q$, the functions $\Phi^{\tilde X,M}$ and $\Phi^{\tilde Y,M}$ are
quasi-equivalent; that is,
\begin{equation*}
  \Phi^{\tilde X,M}(v) \le A\cdot\Phi^{\tilde Y,M}(Bv) + Cv + D
\end{equation*}
for some constants $A$, $B$, $C$, $D$, and \emph{vice versa}. The chain
versions $\Phi^{\tilde X,q}$ and $\Phi^{\tilde Y,q}$ are also quasi-equivalent.
Moreover, if $X$ is a compact Riemannian manifold or compact Lipschitz
neighborhood retract (CLNR) with a triangulation, the profile $\Phi^{\tilde X,M}$
may be interpreted as applying to cellular maps or to Lipschitz maps; the two
interpretations are quasi-equivalent functions, and similarly for $\Phi^{\tilde X,q}$.

As in \cite{jmA99} (with $M=D^q$) and \cite{nB09}, given a finite CW complex or CLNR $X$
where $\tilde X$ is $(q-1)$-connected, one can define the isoperimetric
profiles of $G=\pi_1(X)$ to be those of $\tilde X$; that is,
$\delta_G^M=\Phi^{G,M}=\Phi^{\tilde X,M}$ and $FV_G^q=\Phi^{G,q}=\Phi^{\tilde X,q}$.
The above results ensure that $\Phi^{G,M}$ and $\Phi^{G,q}$ are
well-defined up to quasi-equivalence. It remains to ask whether
$\Phi^{G,M}=\Phi^{G,q}$, or whether $\Phi^{G,M}=\Phi^{G,N}$ for $N$ another
$q$-dimensional connected orientable manifold with nonempty boundary.

Our major result is that $\Phi^{G,M}=\Phi^{G,q}$ almost everywhere
for $q\ge 4$ and $M$ a
$q$-dimensional manifold, and that $\Phi^{G,M}=\Phi^{G,N}$ for $M$ and
$N$ manifolds of dimension $q\ge 3$ with $\bdy M=\bdy N$. More generally,
$\Phi^{X,M}=\Phi^{X,q}$ almost everywhere and $\Phi^{X,M}=\Phi^{X,N}$
under these conditions, provided $X$ is $(q-1)$-connected. In the cellular case, these equalities
are exact; in the metric case, the second is exact, while the first may fail
at discontinuities of the functions (of which there are at most countably
many, since both functions are increasing). This generalizes a theorem
in \cite{nB09}, which states that $\Phi^{G,M^q}\le\Phi^{G,D^q}$ provided $q\ge 4$
and either $\bdy M$ is connected or $\Phi^{G,D^q}$ is superadditive, and which
can be proved by itself by the methods of lemma~\ref{L:singleSimplex2}.

The proof uses, and is inspired by, a classic theorem of Brian White, proved in \cite{bW84}:
Let $X$ be a simply connected Riemannian manifold,
let $M^q$ be an orientable connected manifold with boundary where $q\ge 3$,
and let $f\colon M\to X$ be Lipschitz. Then the Plateau problem for Lipschitz
maps $g\homot f$ rel $\bdy M$ is equivalent to the Plateau problem for
integral currents $T\sim f_\sharp([M])$.

Since we are looking at all functions
which fill a given map $f\colon\bdy M\to X$, and all currents $T$ where
$\bdy T=f_\sharp([M])$, this theorem is not quite sufficient. To finish the proof,
we need for any integral current $T$ a function
$f$ where $f_\sharp([M])$ is very close to $T$.
We use the assumption that $X$ is
$(q-1)$-connected, the Hurewicz maps, and the strong approximation theorem
to prove that such a function must
exist. As in \cite{bW84}, we start with $X$ a CW complex (a mild generalization,
as \cite{bW84} uses simplicial complexes) and use this to prove the metric case.

Section~1 establishes the above argument in detail, allowing
us in section~2 to define the isoperimetric profiles for certain groups $G$
and establish their equality. In a third section, we show that the profiles
for a given group $G$ are finite everywhere, and computably bounded
provided $G$ has solvable word problem.

Most of the concepts and conventions used herein were introduced in
\cite{clG08}, and the acknowledgements there apply here as well.

\section{The map \emph{vs.} current equivalence}

Let $X$ be a $(q-1)$-connected CW complex where $q\ge 3$,
and let $M$ be a $q$-dimensional compact orientable manifold
with $\bdy M\ne\emptyset$.
We want to show that $\Phicell^{X,M} \le \Phich^{X,q}$, with
equality in dimensions $q\ge 4$.
To do this, we must show that
chains in dimensions $q$ and $q-1$ can be represented by
functions $M^{q-1}\to X^{(q-1)}$ and $M^q\to X^q$, where $M$ is
some compact manifold, possibly with boundary. The easiest cases
are where $M=D^q$ or $M=S^{q-1}$; here we can interpret ``volume''
to mean ``word length in $\pi_q(X^{(q)},X^{(q-1)})$'' or the same
one dimension down. Since $X$ is highly connected, the Hurewicz
maps are length-preserving isomorphisms, and this is sufficient.
For more general $M$, we can triangulate $M$ and put all the volume
on a single cell.

The metric case is similar, with one wrinkle. If $T$ is a polygonal
current in $X$ of dimension $q$, then $T$ can be represented by a
Lipschitz map $f\colon M\to X$ with the same volume (by the above logic
applied to some highly connected simplicial complex which supports $T$).
However, an integral current $T$ with volume $v$ can only be approximated by
a polygonal current, with volume at most $v+\epsilon$. Thus, if $\Phimet^{X,M}$
has a jump discontinuity at $v$, then $\Phicur^{X,q}(v)$ lies between
$\Phimet^{X,M}(v)$ and the right limit. As $\Phimet^{X,M}$ is an increasing function,
this implies that $\Phimet^{X,M}=\Phicur^{X,q}$ almost everywhere,
and certainly the functions are quasi-equivalent.

To begin, consider the diagram
\begin{equation*}\label{E:Hurewicz}
  \xymatrix{
    \pi_q(X^{(q)},X^{(q-1)}) \ar[r]^\varphi_\sim \ar[d]_\bdy &
    H_q(X^{(q)},X^{(q-1)}) \ar[d]^\bdy \\
    \pi_{q-1}(X^{(q-1)}) \ar[r]^\varphi_\sim \ar@{^{(}->}[d]_{j_*} &
    H_{q-1}(X^{(q-1)}) \ar@{^{(}->}[d]^{j_*} \\
    \pi_{q-1}(X^{(q-1)},X^{(q-2)}) \ar[r]^\varphi &
    H_{q-1}(X^{(q-1)},X^{(q-2)})
  }
\end{equation*}
where $j\colon (X^{(q-1)},*)\to(X^{(q-1)},X^{(q-2)})$ is the
inclusion, each $\varphi$ is a Hurewicz homomorphism, and each $\bdy$
is the usual boundary map.  By both parts of \cite[theorem
7.4.3]{ehS66}, this diagram commutes.  By the Hurewicz isomorphism
theorem, all of the maps $\varphi$ are isomorphisms if $q\ge 4$, and
the top and bottom maps preserve word length.  If $q=3$, then the top
two $\varphi$ are still isomorphisms, the top $\varphi$ still
preserves length, and the last $\varphi$ is a length-reducing
epimorphism.  The $j_*$ on the right is an injection, by the long
exact sequence for $(X^{(q-1)},X^{(q-2)})$ (certainly
$H_{q-1}(X^{(q-2)})=0$), so the $j_*$ on the left is an injection as
well for $q\ge 4$.  (The left $j_*$ is also an injection if $q=3$ by
the exact sequence of homotopy groups, since $\pi_2(X^{(1)})=0$.)
Note that the composition $j_*\comp\bdy$ on the right is our previous
notion of $\bdy$ on chains; since $j_*$ is 1-1, we might ignore this
distinction.

\begin{lemma}\label{L:singleSimplex2}
  Let $q\ge 3$, let $\tau$ be a triangulation of $M$, and let $X$ be a
  $(q-1)$-connected CW complex.  Let $f\colon\bdy M\to X$ be
  quasi-cellular, and let $T$ be a $q$-chain of $X$ where $\bdy
  T=f_\sharp([\bdy M])$.  Then there exists a $\tau$-cellular map
  $g\colon M\to X$ where $\bdy g=f$, $g_\sharp([M])=T$, and $\Vol_\tau
  g = \|T\|$.
\end{lemma}

\begin{proof}
  Recall that $\Vol_\tau g$ is the sum, over all $q$-cells $\Delta$
  of the triangulation $\tau$, of the word length of $[g\restriction\Delta]$
  in $\pi_q(X^{(q)},X^{(q-1)})$.
  In fact, we will construct $g$ so that every $q$-cell except
  possibly one is sent to $X^{(q-1)}$. Let $G=(V,E)$ be the
  undirected graph where $V$ is the set of $q$-cells of $\tau$,
  and where $\{v,w\}\in E$ iff $v$ and $w$ share a $(q-1)$-face.
  $G$ is connected, so let $R$ first be a spanning tree for $G$.
  We define a procedure which will define $g$ cell by cell on $M$,
  while removing cells from $R$.
  At any given stage, $g$ will be defined
  only on $\bdy M$ and on those cells which are no longer in $R$.  In
  particular, if $\{v,w\}\in R$ at any stage, then $g$ will not be
  defined on the interior of $v\cap w$, as this cell is neither a
  subset of $\bdy M$ nor a face of any $q$-cell where $g$ is already
  defined.

  Proceed as follows: While $R$ contains more than one vertex, let $v$
  be a leaf of $R$ and let $w$ be the unique $q$-cell where $\{v,w\}\in R$.
  Let $D=(\bdy v)\setminus (v\cap w)^\circ$, which is homeomorphic to
  $D^{q-1}$.  There is
  some subcomplex of $D$ on which $g$ is already defined (possibly
  empty).  Since $X^{(q-1)}$ is $(q-2)$-connected, we may extend $g$
  in some way to all of $D$.  Finally, there is a retraction $r\colon
  v\to D$; define $g$ on $v$ as $g\comp r$.  (Note that
  $g_\sharp([v])=g_\sharp(r_\sharp([v]))=0$.)  Now that $g$ is defined
  on $v$, remove the vertex $v$ and the edge $\{v,w\}$ from $R$.  The
  remaining cells still form a connected, \emph{etc.}~manifold, at
  least in the PL category, so we may repeat.

  When this procedure is finished, let $\Delta$ be the unique $q$-cell
  remaining in $R$, so that $g$ is defined on
  $M\setminus\Delta^\circ$.  Since $g_\sharp([v])=0$ for every other
  $q$-cell $v$, we have $g_\sharp([\bdy\Delta])=\bdy T$.  Choose
  $h\colon(D^q,S^{q-1})\to(X^{(q)},X^{(q-1)})$ where $\varphi([h])=T$.
  Then there is a homotopy $H\colon \bdy h\simeq f$.  Attach $H$ to
  $h$ to obtain $g\restriction\Delta$.
\end{proof}

\begin{theorem}\label{T:trEqch}
  Let $q\ge 4$ and $\bdy M\ne\emptyset$, and let $X$ be a
  $(q-1)$-connected CW complex.  Then $\Phich^{X,q}=\Phicell^{X,M}$.  If
  $q=3$, then $\Phicell^{X,M}\le\Phich^{X,q}$.
\end{theorem}

\begin{proof}
  Let $q\ge 3$, let $f\colon \bdy M\to X^{(q-1)}$ be a
  quasi-cellular map, and suppose $\Voltr f\le n$.  There is a
  triangulation $\tau$ of $M$ where $\Vol_{\bdy\tau} f=\Voltr f$.
  Because $X^{(q)}$ is $(q-1)$-connected, there is some extension of
  $f$ to a $\tau$-cellular map $f'\colon M\to X$.  Thus the
  $(q-1)$-chain $S=f_\sharp([\bdy M])$ is the boundary of the
  $q$-chain $f'_\sharp([M])$.  Let $T$ be a $q$-chain of least
  possible volume where $\bdy T=S$.  By lemma \ref{L:singleSimplex2},
  there is a $\tau$-cellular map $g\colon M\to X$ where $\bdy g=f$ and
  $\Vol_\tau g = \|T\| \le \Phich^{X,q}(n)$.  Generalizing over all
  $f$, $\tau$, and $n$, $\Phicell^{X,M} \le \Phich^{X,q}$.

  To see the reverse inequality for $q\ge 4$, choose a
  triangulation $\tau$ on $M$ and a $q$-cell $\Delta$ with at least
  one face on the boundary on $M$; call this face $\delta$.  Let $s$
  be the map on $M$ which collapses all points of $M$, except those in
  $\Delta^\circ$ or $\delta^\circ$, to a single point $*$; $s$ is a
  map from $M$ to $D^q$ which is a diffeomorphism $\Delta^\circ\cong
  (D^q)^\circ$ and $\delta^\circ\cong (S^{q-1}\setminus\{*\})$.

  Given a $(q-1)$-boundary $S=\bdy T$ where $\|S\|\le n$, let
  $f'\colon S^{q-1}\to X^{(q-1)}$ where $\varphi([f'])=S$ and
  $f'(*)\in X^{(0)}$ (this is where we use $q\ge 4$).  Let $f=f'\comp s\colon \bdy M\to X^{(q-1)}$,
  so that $\Vol_\tau f = \|S\|$.  $f'$
  has a filling $g'$, by the Hurewicz isomorphism, so $f$ has a
  filling $g'\comp s$.  Choose a filling $g$ for $f$ so that $\Voltr
  g$ is minimum; then $T=g_\sharp([M])$ fills $S$ and
  \begin{equation*}
    \FVch(S) \le \|T\| \le \Voltr g = \FVcell f \le \Phicell^{X,M}(n).
  \end{equation*}
  Generalizing over all $S$ and $n$, $\Phich^{X,q}\le\Phicell^{X,M}$.
\end{proof}

\begin{corollary}
  For $q\ge 4$, $X$ as above, and $M$, $N$ manifolds with nonempty
  boundary, $\Phicell^{X,M}=\Phicell^{X,N}$.  This still holds true for
  $q=3$ provided $\bdy M\cong\bdy N$.
\end{corollary}

\begin{proof}
  The case $q\ge4$ follows directly from theorem \ref{T:trEqch};
  $\Phicell^{X,M}=\Phich^{X,q}=\Phicell^{X,N}$.  For $q=3$, refer to the
  first paragraph of the above proof; by this reasoning, every
  quasi-cellular function $f\colon (\bdy M=\bdy N)\to X^{(q-1)}$ has a
  filling on both $M$ and $N$, and moreover $\FVcell^M f = \FVch
  f_\sharp([\bdy M]) = \FVcell^N f$.  Thus $\Phicell^{X,M} =
  \Phicell^{X,N}$.
\end{proof}

Now we address the metric case. Let $X$ be an LLNR;
that is, let $X\subseteq U \subseteq \Rn^N$ where $U$ is open in $\Rn^N$,
and let $r\colon U\to X$ be a locally Lipschitz retraction.

\begin{theorem}\label{T:curRepresent}
  Let $q\ge3$, $\bdy M\ne\emptyset$, $X$ a $(q-1)$-connected LLNR, $P$
  a finite simplicial complex of dimension $q$, and $\psi\colon P\to
  X$ a 1-1 Lipschitz map.  Let $T=\psi_\sharp([P])$.  Then there
  exists $f\colon M\to X$ where $f_\sharp([M])=T$ and $\Volmet f=\Volcur(T)$.
  If $q\ge 4$, one may choose $f$ so that $\Volmet\bdy f=\Volcur(\bdy T)$.
\end{theorem}

\begin{proof}
  Construct $Q\supset P$ so that $Q$ is $q$-dimensional and
  $(q-1)$-connected.  (For example, find a high-dimensional simplex
  which contains $P$ as a subcomplex and take its $q$-skeleton.) 
  Extend $\psi$ to a continuous map $\psi'\colon Q\to X$; this is possible
  because $X$ is $(q-1)$-connected. Mollify $\psi'$ outside $P$ to
  obtain a Lipschitz map $\psi''\colon Q\to U$, and let $\psi'''=r\comp\psi''$.
  Then $\psi'''$ is a Lipschitz map $Q\to X$ which extends $\psi$.
  For convenience, refer to $\psi'''$ as $\psi$.

  As in the proof of theorem \ref{T:trEqch}, choose a triangulation
  $\tau$ of $M$, a cell $\Delta$ intersecting $\bdy M$, and a
  collapsing map $s\colon(M,\bdy M)\to(D^q,S^{q-1})$ which is a
  diffeomorphism on $\Delta^\circ$ and $(\Delta\cap\bdy M)^\circ$.
  Choose $g\colon(D^q,S^{q-1})\to(Q^{(q)},Q^{(q-1)})$ where
  $\varphi([g])=[P]$.  We may assume
  that $g$ covers each point in the interior of a $q$-cell of $P$
  exactly once, that it does so smoothly, and that $g[D^q]$ does
  not intersect the interior of
  any other $q$-cell of $Q$; and similarly for $\bdy g$
  provided $q\ge 4$.  Let $f=\psi\comp g\comp s$.  Then
  \begin{equation*}
    f_\sharp([M]) = \psi_\sharp(g_\sharp([D^q])) = \psi_\sharp([P]) = T.
  \end{equation*}
  Over every cell of $\tau$ except for $\Delta$, the volume form of
  $f^*(ds^2)$ disappears, and we calculate
  \begin{equation*}
    \Volmet f = \Volmet (\psi\comp g) = \int_P \psi^*(ds^2) = \Volcur(T),
  \end{equation*}
  the last because $\psi$ is 1-1.  A similar calculation shows that
  $\Volmet \bdy f = \Volcur(\bdy T)$ for $q\ge 4$.
\end{proof}

A similar result holds for $q\ge 3$ and $\bdy M=\emptyset$, provided
$\bdy T=0$.  (The only way this can happen is if $\bdy P=\emptyset$ as
well, which by exactness implies $P=j_*(P')$ for some $P'\in
H_{q-1}(Q^{(q-1)})\cong\pi_{q-1}(Q^{(q-1)})$.)

\begin{theorem}\label{T:curEqTrAE}
  Let $q\ge 4$ and let $M$ and $X$ be as in theorem
  \ref{T:curRepresent}.  Then
  \begin{equation*}
    \Phimet^{X,M}\le\Phicur^{X,q}\le\bar{\Phimet^{X,M}}.
  \end{equation*}
  In particular, $\Phimet^{X,M}=\Phicur^{X,q}$ almost everywhere.  If
  $q=3$, then $\Phimet^{X,M}\le\Phicur^{X,q}$.
\end{theorem}

Recall that $\bar f$ is the upper envelope of $f$, where
$f$ is any map from a topological space to $\Rn$.

\begin{proof}
  \cite[theorem 3]{bW84} tells us that, given a function $f\colon M\to
  X$,
  \begin{equation*}
    \inf\{\,\Volmet g : g\simeq f\textrm{ rel $\bdy M$}\,\}
    = \inf\{\,\Volcur(T) : T - f_\sharp([M])\in\mathcal{B}_q(X)\,\}.
  \end{equation*}
  Let $q\ge 3$, $r\ge0$, and $f\in C^{0,1}(\bdy M,X)$, $\Volmet f\le r$.
  Since $X$ is $(q-1)$-connected, $f$ is the boundary of some
  continuous map $h\colon M\to X$, and $h$ can be mollified and retracted
  to a Lipschitz map $M\to X$.  Thus $S=f_\sharp([\bdy M])$ is the
  boundary of $h_\sharp([M])$.
  
  Let $R\in\mathbf{I}_q(X)$ where
  $\Volcur(R)<\FVcur S + \epsilon$.  $R-h_\sharp([M])$ is a $q$-cycle,
  which is a $q$-boundary up to an element of $H_q(X)$ (singular
  homology).  But $H_q(X)\cong\pi_q(X)$, so by modifying $h$ in a
  disk, we may assume $R-h_\sharp([M])$ is a boundary.  Thus
  \begin{align*}
    \FVmet f &\le \inf\{\,\Volmet g : g\simeq h\textrm{ rel $\bdy M$}\,\}
    = \inf\{\,\Volcur(T): T-h_\sharp([M])\in\mathcal{B}_q(X)\,\} \\
    &\le \Volcur(R) < \Phicur^{X,q}(r) + \epsilon.
  \end{align*}
  Taking $\epsilon\to0$ and generalizing over all $f$ and $r$, we have
  $\Phimet^{X,M}\le\Phicur^{X,q}$.

  Conversely, for $q\ge 4$ let $S\in\mathbf{I}_{q-1}(X)$ where $\bdy
  S=0$ and $\Volcur(S)\le r$.  By the strong approximation theorem
  \cite[lemma 4.2.19]{hF69}, $S$ is homologous by a current $R$ where
  $\Volcur(R)<\epsilon$ to a polyhedral current $S'=\psi_\sharp([P])$,
  $\Volcur(S')<\Volcur(S)+\epsilon$.  Note that $\bdy S'=0$ as well.  By
  theorem \ref{T:curRepresent}, choose $f\colon\bdy M\to X$ where
  $f_\sharp(\bdy M]) = S'$ and $\Volmet f=\Volcur(S')$.  $f$ extends to
  some $g\colon M\to X$; choose $g$ where $\Volmet g <
  \FVmet(f)+\epsilon$.  Then
  \begin{equation*}
    \FVcur S < \Volcur(g_\sharp([M])) + \epsilon
    < \FVmet(f) + 2\epsilon \le \Phimet^{X,M}(r+\epsilon) + 2\epsilon;
  \end{equation*}
  taking $\epsilon\to0$, $\FVcur S\le \bar{\Phimet^{X,M}}(r)$ (since
  $\Phimet^{X,M}$ is increasing).
  Generalizing over $S$ and $r$, we have
  $\Phicur^{X,q}\le\bar{\Phimet^{X,M}}$.

  A function can differ from its upper envelope only at points of
  discontinuity; since $\Phimet^{X,M}$ is increasing, there are at most
  countably many of these (see \cite{hlR88}), so
  $\Phimet^{X,M}=\Phicur^{X,q}$ almost everywhere.
\end{proof}

\begin{corollary}\label{C:sameManifold}
  For $q\ge 4$, $X$ as in theorem \ref{T:curRepresent}, and $M$ and
  $N$ manifolds with nonempty boundary, $\Phimet^{X,M}=\Phimet^{X,N}$
  almost everywhere.  Provided $\bdy M\cong\bdy N$, the conclusion
  holds for $q\ge 3$ and with exact equality.
\end{corollary}

\begin{proof}
  The first follows from theorem \ref{T:curEqTrAE}.  For the
  second, follow the reasoning in the first paragraph of the
  preceding proof.  Every $f\in C^{0,1}(\bdy M,X)$ has a filling, and
  $\FVmet^M f = \FVcur f_\sharp([\bdy M]) = \FVmet^N f$.  Therefore
  $\Phimet^{X,M}=\Phimet^{X,N}$.
\end{proof}

The case $q=3$ deserves some attention. For simplicity,
assume that we only consider $M$ where $\bdy M$ is connected;
thus $\bdy M=\Sigma_g$ is the surface of genus $g$ for some $g\ge 0$.
By Corollary~\ref{C:sameManifold}, only $\bdy M$ is relevant,
so assume $M=\Gamma_g$ is the solid torus of genus $g$.
If $g\le h$, then every map $\Gamma_g\to X$ can be composed with
a collapsing map $\Gamma_h\to\Gamma_g$ to produce a new map with equal
volume and filling volume. Thus
$\Phicell^{X,\Gamma_g}\le\Phicell^{X,\Gamma_h}$
and $\Phimet^{X,\Gamma_g}\le\Phicell^{X,\Gamma_h}$ everywhere.
As every 2-chain or polygonal 2-current $T$ can be represented
by a map $f\colon\Sigma_g\to X$ for some $g$, we have
$\Phich^{X,3} = \lim_{g \uparrow 0} \Phicell^{X,\Gamma_g}$ everywhere
and $\Phicur^{X,3} = \lim_{g \uparrow 0} \Phimet^{X,\Gamma_g}$
almost everywhere.

As we will note later, there are spaces $X$ for which
$\Phicell^{X,D^3}\ne\Phicell^{X,\Gamma_1}$. Such a separation
between $\Phicell^{X,\Gamma_g}$ and $\Phicell^{X,\Gamma_h}$
for $1\le g<h$, or between $\Phicell^{X,\Gamma_g}$ and $\Phich^{X,3}$,
is not yet known.

\section{Applications to geometric group theory}

We say that a group $G$ is of type $F_q$ if there is a CW
complex $X=K(G,1)$ where $X^{(q)}$ is finite. Equivalently,
$G$ is of type $F_q$ iff there is a finite CW-complex $Y$
where $\tilde Y$ is $(q-1)$-connected and $\pi_1(Y)=G$.
$Y$ may be taken as $q$-dimensional; alternatively,
$Y$ may be a compact manifold.

For example, every group is of type $F_0$. A group is of
type $F_1$ iff it is finitely generated, and of type $F_2$ iff
it is finitely presented. Every type $F_q$ is a strict subtype of
$F_{q-1}$; this follows in the case $q=3$ by the results in
\cite{jS63}.

In this case, where $X=K(G,1)$ has finite $q$-skeleton,
we may take the functions $\Phich^{\tilde X,q}$ and
$\Phicell^{\tilde X,M}$ (for $M$ a $q$-manifold) as invariants
of $G$. For example, $\Phicell^{\tilde X,D^2}$ is the classical
Dehn function, while $\Phicell^{\tilde X,D^q}$ is the higher-order
Dehn function $\delta_{q-1}$ studied in \cite{jmA99}.
As in these cases (and for similar reason), if we change $X$
to a different $K(G,1)$, we obtain a quasi-equivalent function.
By the results of the last section, many of these functions
coincide.

\begin{lemma}\label{L:groupqEquiv}
  Let $q\ge 2$, and let $G$ be a group of type $F_q$.
  Let $X$ and $Y$ be CW complexes which are $K(G,1)$'s
  and where $X^{(q)}$ and $Y^{(q)}$ are both finite.  Then
  $\Phich^{\tilde X,q}\qequiv\Phich^{\tilde Y,q}$.  Also
  $\Phicell^{\tilde X,M}\qequiv\Phicell^{\tilde Y,M}$ for $M$
  $q$-dimensional and $\bdy M\ne\emptyset$.
\end{lemma}

\begin{proof}
  Let $\varphi\colon \pi_1(X,*)\to\pi_1(Y,*)$ be an isomorphism.  By
  \cite[Theorem 1B.9]{aH02}, there is a continuous function $f\colon
  X\to Y$ where $f_*=\varphi$.  The conclusions follow from 
  \cite[theorem 2]{clG08}.
\end{proof}

\begin{definition}
  Let $q\ge 2$. Let $G$ be a group of type $F_q$, and let $X$ be a $K(G,1)$ where
  $X^{(q)}$ is finite. The
  \emph{chain isoperimetric profile} of $G$ in dimension $q$ is that of $\tilde X$;
  \begin{equation*}
    \Phich^{G,q} := \Phich^{\tilde X,q}.
  \end{equation*}
  Likewise, for $M$ a $q$-manifold with $\bdy M\ne\emptyset$,  the
  \emph{cellular isoperimetric profile} of $G$ for $M$ is that of $\tilde X$;
  \begin{equation*}
    \Phicell^{G,M} := \Phicell^{\tilde X,M}.
  \end{equation*}
\end{definition}

For example, if $G$ is finitely presented, the function $\Phicell^{G,D^2}$ is
defined, and is in fact the usual Dehn function of $G$. As one would expect,
this definition is an abuse of language, since $\Phich^{G,q}$ and $\Phicell^{G,M}$
are defined only up to quasi-equivalence.
It still makes sense to ask whether these functions are linear, or polynomial,
or exponential, or computably bounded, or everywhere finite.

\begin{theorem}
 Let $G$ is a group of type $F_q$, where $q\ge 4$.  Then for
  any $q$-manifold $M$ with $\bdy M\ne\emptyset$,
  $\Phicell^{G,M}\qequiv\Phich^{G,q}$.  If $q=3$, then
  $\Phicell^{G,M}\qbound\Phich^{G,q}$; and if $M$ and $N$
  are two 3-manifolds with $\bdy M=\bdy N\ne\emptyset$, then
  $\Phicell^{G,M}\qequiv\Phicell^{G,N}$.
\end{theorem}

\begin{proof}
  Apply theorem~\ref{T:trEqch} and corollary~\ref{C:sameManifold} to any
  $K(G,1)$ with finite $q$-skeleton.
\end{proof}

\begin{theorem}
  Let $X$ be a closed Riemannian manifold where $\pi_1(X)=G$ and
  $\tilde X$ is $(q-1)$-connected.  Then $\Phicur^{\tilde
    X,q}\qequiv\Phich^{G,q}$, and $\Phimet^{\tilde
    X,M}\qequiv\Phicell^{G,M}$ for any $M^q$.
\end{theorem}

\begin{proof}
  Choose a triangulation $\tau$ on $X$.  There is a $K(G,1)$ whose
  $q$-skeleton is $\tau^{(q)}$; one adds $(q+1)$-cells to kill
  $\pi_q$, \emph{etc.}  Thus
  \begin{equation*}
    \Phicur^{\tilde
    X,q}\qequiv\Phich^{\tilde\tau,q}\qequiv\Phich^{G,q}\textrm{ and }
    \Phimet^{\tilde X,M}\qequiv\Phicell^{\tilde\tau,M}\qequiv\Phicell^{G,M}.
    \qedhere
  \end{equation*}
\end{proof}

\section{Finiteness and computability}\label{S:finComp}

One of the classic results on Dehn functions of finitely presented
groups $G$ is that $\delta_G$ is recursive, or even subrecursive,
iff the word problem on $G$ is solvable. The proof from right to
left is fairly simple: Fix $X=K(G,1)$ with finite 2-skeleton. Up to
translation by $G$, there are only finitely many loops in $\tilde X$
whose length is at most $n$, for any fixed $n$. Of these, use the word
problem solution to determine which are contractible, then search
all the disk maps into $\tilde X$ (up to translation)
until one is found for each contractible loop.

The situation for generalized
Dehn functions is not so clear-cut. For example, a result by Papasoglu
(see \cite{pP00}) states that the ``second Dehn function'' $\Phicell^{G,D^3}$
is subrecursive for \emph{all} groups $G$ of type $F_3$, while an unpublished
result by Young (see \cite{rY08}) shows that there is a group $G$, which
is of type $F_q$ for all $q$, where
$\Phicell^{G^2,\Gamma_1}$ is not subrecursive, nor is $\Phich^{G^{q-1},q}$ for
any $q\ge 4$. In fact, the second Dehn function is effectively the only generalized
Dehn function which is subrecursive for all qualifying $G$.

Nevertheless, results in one direction is possible.

\begin{theorem}\label{T:boundedPhiCh}
  Let $q\ge 2$ and let $G$ be a group of type $F_q$. Then
  $\Phich^{G,q}(n)$ is finite for all $n$. If $G$ has solvable word
  problem, then $\Phich^{G,q}$ is recursive.
\end{theorem}

The proof is fairly similar to that for the classical Dehn function, with the wrinkle
that the class of chains with volume at most $n$ is not generally finite, even
up to translation by $G$. We observe, however, that each chain can be
decomposed into connected components, and the class of small
connected chains can be exhausted. The equivalent of $\Phich^{G,q}$
where only connected boundaries are considered is therefore computable.
If the components of a cycle are far enough apart, then the most efficient
way to fill the cycle is to fill each component separately. This allows us to
compute $\Phich^{G,q}$ in terms of the specialized version.

\begin{theorem}\label{T:boundedPhi}
  Let $q\ge 2$, and
  let $G$ be a group of type $F_q$, and
  let $M$ be a $q$-dimensional manifold with nonempty boundary. Then
  $\Phicell^{G,M}(n)$ is finite for all $n$. If $G$ has solvable word problem, then
  $\Phicell^{G,M}$ is subrecursive.
\end{theorem}

For $q\ge 3$, the [recursive] bound on $\Phicell^{G,M}$ is precisely $\Phich^{G,q}$;
for $q=2$, the bound may be computed in terms of the classical Dehn function.

We need some preliminaries. Let $X$ be a CW-complex.
For any $t$ and for any $t$-chains $A$, $B$, we say $B$
is a \emph{subchain} of $A$ if $\|A\| = \|B\| + \|A-B\|$. Equivalently,
for ever $t$-cell $\sigma$ of $X$, let $n_{A,\sigma}$ and $n_{B,\sigma}$
be the coefficients of $\sigma$ in the expansions of $A$ and $B$ respectively.
Then $B$ is a subchain of $A$ iff $0\le n_{B,\sigma}\le n_{A,\sigma}$
or $n_{A,\sigma}\le n_{B,\sigma}\le 0$ for all $\sigma$. A given chain $A$
has finitely many subchains.

We say $B$ is a \emph{component} of $A$ if $B$ is a subchain of $A$
and $\bdy B$ is a subchain of $\bdy A$, and that a $t$-chain $A$ is
\emph{connected} if its only components are 0 and $A$. By induction
on $\|A\|$, every chain $A$ can be expressed as a sum of connected
components $A=B_1+\dots+B_n$. Given such an expansion, we have
$\|A\|=\|B_1\|+\dots+\|B_n\|$ and
$\|\bdy A\| = \|\bdy B_1\| + \dots + \|\bdy B_n\|$. Thus, if $\bdy A=0$,
then $\bdy B_i=0$ for all $i$.

\begin{proof}[Proof of theorem~\ref{T:boundedPhiCh}]
  Let $X$ be a $K(G,1)$ with finite $q$-skeleton, and let $\tilde X$
  be its universal cover. Assume $G$ is finite, so that $\tilde X^{(q)}$ is finite.
  For $0\le t\le q$ and $v\ge 0$, we can determine the set of all
  elements $T\in C_t(\tilde X)$ where $\|T\|\le v$. Also we can
  determine whether a given $T\in C_t(\tilde X)$ is a cycle, hence a
  boundary. Thus we have the following algorithm for $\Phich^{G,q}(n)$:
  determine all the $(q-1)$-chains $S$ which are boundaries; for each
  $N$ starting with 0 and every $q$-chain $T$ where
  $\|T\|=N$, flag $\bdy T$; stop when every $S$ is flagged and output $N$.
  
  Now assume $G$ is infinite. Each cell of $X$ is covered by a collection of
  cells of $\tilde X$ in 1-1 correspondence with the elements of $G$,
  and the collection is invariant under the natural $G$-action. For each $t\le q$,
  let $\Sigma_t$ be a
  set containing exactly one $t$-cell of $\tilde X$ which covers any given
  $t$-cell of $X$; thus $\Sigma_t$ is finite.
  
  Up to $G$-action, there are finitely many connected $t$-chains of any fixed
  volume $n$. For $n=0$, this is obvious. For $n\ge 1$, suppose $A$ is a
  connected $t$-chain and $\|A\|=n$. We generate chains $B_1$, $B_2$, \dots,
  $B_n=A$, each a subchain of the next, where $\|B_k\|=k$, as follows: After translation by
  some element of $G$, $A$ must contain $\pm\sigma$ where $\sigma\in\Sigma_t$.
  Let $B_1=\pm\sigma$. Suppose $B_k$ has been chosen for some $k<n$.
  Since $A$ is connected, $B_k$ is not a component of $A$, so there is some
  subchain $C$ of $A-B_k$ where $\|C\|=1$ and $\|\bdy(B_k+C)\|<\|\bdy B_k+\bdy C\|$.
  Let $B_{k+1}=B_k+C$. There are only finitely many possibilities for $B_1$; and
  for any chain $B$, there are only finitely many cells which share boundary with $B$,
  so for any choice of $B_k$ there are finitely many choices for $B_{k+1}$. Thus there
  are finitely many possible $B_n=A$, which is what we wanted.
  
  Let $\Psi(n)=\max \FVch(A)$, where the maximum is taken over
  \emph{connected} $(q-1)$-cycles $A$ with $\|A\|\le n$. Let
  \begin{equation*}
    \Phi(n) = \max_\textrm{partitions $P$ of $n$} \sum_{k\in P} \Psi(k),
  \end{equation*}
  where $P$ is interpreted as a multiset. $\Psi(n)$ and $\Phi(n)$ are
  finite for all $n$.
  
  We claim that $\Phich^{\tilde X,q}(n)=\Phi(n)$ for all $n$.
  To see $\Phich^{\tilde  X,q}(n)\le\Phi(n)$, let $A$ be a
  $(q-1)$-cycle with $\|A\|\le n$. Let $A=B_1+\dots+B_m$ be a
  sum of connected components. Each $B_i$ is a cycle; let
  $C_i\in C_q(\tilde X)$ with minimum volume where $\bdy C_i=B_i$.
  The multiset $\{\|B_1\|,\dots,\|B_m\|\}$ is a partition of $\|A\|$, so
  \begin{equation*}
    \FVch(A) \le \sum_{i=1}^m \|C_i\| = \sum_{i=1}^m \FVch(B_i)
    \le \sum_{i=1}^m \Psi(\|B_i\|) \le \Phi(\|A\|)\le\Phi(n).
  \end{equation*}
  Take the supremum over all $A$ to see $\Phich^{\tilde X,q}(n)\le\Phi(n)$.
  In particular, $\Phich^{\tilde X,q}(n)$ is finite for all $n$.
  
  To see $\Phi(n)\le\Phich^{\tilde X,q}(n)$, let $B_1$, \dots, $B_m$ be
  a finite list of connected $(q-1)$-cycles where $\sum_i \|B_i\|\le n$.
  Choose translates $B'_i=g_iB_i$ as follows: Let $g_1=e$ and
  $B'_1=B_1$. For $i=j+1$, note that there are finitely many connected
  $q$-chains $C$ where $\bdy C$ shares a cell with any of $B'_1$, \dots, $B'_j$
  and where $\|C\|\le\Phich^{\tilde X,q}(n)$. Since $G$ is infinite, there must
  be some $h\in G$ where $h=B_i$ does not share a cell with
  $\bdy C$ for any such $C$. Take $g_i$ to be some such $h$.
  Note that $\|B'_i\|=\|B_i\|$, $\bdy B'_i=0$, and $\FVch(B'_i)=\FVch(B_i)$
  for all $i$.
  
  Let $A=B'_1+\dots+B'_m$, so that $\bdy A=0$ and $\|A\|\le n$.
  Let $C$ be a $q$-chain where $\bdy C=A$ and $\|C\|=\FVch(A)$.
  If $C'\ne0$ is a connected component of $C$, then $\bdy C'$ shares cells
  with $B'_i$ for a unique $i$, since $\|C'\|\le\|C\|\le\Phich^{\tilde X,q}(n)$
  and we chose the $B'_i$ so that no such connected $q$-chain exists
  with boundary cells from two distinct $B'_i$. Conversely, any boundary
  cells of $\bdy C'$ must come from $\bdy C=A=\sum_i B'_i$, and if $\bdy C'=0$,
  then $C-C'$ is a filling cycle for $A$ of smaller volume, which is impossible.
  Thus, given a sum of connected components $C=C_1+\dots+C_N$, each
  $C_j$ is associated with a unique $B'_i$, and $C$ is therefore a sum over $i$
  of fillings for $B'_i$. Hence
  \begin{equation*}
    \sum_{i=1}^m \FVch(B'_i) = \FVch(A) \le \Phich^{\tilde X,q}(n).
  \end{equation*}
  Taking the supremum over all finite lists $B_1$, \dots, $B_m$, we have
  $\Phi(n) \le \Phich^{\tilde X,q}(n)$, as desired.
  
  Finally, suppose $G$ is infinite with solvable word problem (let $Y$ be
  a finite set of generators). We first
  show that $\Psi$ is computable. If $w\in F(Y)$ is a word and
  $\sigma\in\Sigma_t$ for $t\le q$, then $(w,\sigma)$ represents the $t$-cell
  $\bar w\sigma$. Every $t$-cell is represented by some pair, and $(w,\sigma)$
  and $(w',\sigma')$ represent the same cell iff $\sigma=\sigma'$ and
  $w^{-1}w'$ represents the trivial element of $G$. Chains $A\in C_t(\tilde X)$
  can be represented as linear combinations of pairs. From a representation
  for $A$, one can calculate $\|A\|$; also, one can calculate a representation
  for $\bdy A$. Also, one can determine all the subchains of $A$, and one can
  decide whether a chain $A$ is connected: for each subchain $B$ where $B\ne 0$
  and $B\ne A$, determine whether $\bdy B$ is a subchain of $\bdy A$. $A$
  is connected iff this is never the case.
  
  The proof that $\Psi$ is computable is similar to the case with $G$ finite. Generate
  representations of all of the $(q-1)$-chains $A$ with $\|A\|\le n$,
  up to $G$-action, by starting with cells in $\Sigma_{q-1}$ and adding cells
  which share boundary components; there are finitely many such chains.
  Of these chains, take the subset consisting of connected cycles. For each
  connected cycle $A$, generate the $q$-chains $B$ where $\bdy B$ and $A$
  share at least one cell, in order of increasing volume, stopping when a chain
  $B$ with $\bdy B=A$ is found; remember the volume of $B$. The maximum
  such volume, across all $A$, is $\Psi(n)$.
  
  Given an algorithm for computing $\Psi$, computing $\Phi=\Phi^{\tilde X,q}$
  is easy.
\end{proof}

\begin{proof}[Proof of theorem~\ref{T:boundedPhi}]
  For $q\ge 3$ the theorem is an easy consequence of theorems~\ref{T:trEqch}
  and~\ref{T:boundedPhiCh}. Now consider $q=2$. As noted previously, the
  classical Dehn function of $G$, here $\Phicell^{G,D^2}$, is finite for all inputs,
  and that it is computably bounded (even computable) iff $G$ has
  solvable word problem. For arbitrary $M^2$, we have $\bdy M$
  a disjoint union of $k$ copies of $S^1$, which implies
  \begin{equation*}
    \Phicell^{G,M}(n) \le \max_{n_1+\dots+n_k=n}\sum_{i=1}^k \Phicell^{G,D^2}(n_i).
  \end{equation*}
  Similarly, if $A$ is a 1-cycle, it can be represented as $A_1+\dots+A_k$
  where each $A_i$ is the image of a copy of $S^1$; hence
  \begin{equation*}
    \Phich^{G,2}(n) \le
      \max_\textrm{$P$ a partition of $n$}\sum_{k\in P}\Phicell^{G,D^2}(k).
    \qedhere
  \end{equation*}
\end{proof}

\section{Further questions}
The relationship between $\Phich^{G,2}$ and the various $\Phicell^{G,M}$
for $\dim M=2$ is not yet clear. Also we do not yet know whether $\Phi^{G,M}$ is
always recursive for $\dim M\le 3$ and $G$ decidable, rather than subrecursive. Finally,
there is no obvious converse to the results of section~\ref{S:finComp};
that is, if $G$ does not have
decidable word problem, it is unclear what restrictions this places on
$\Phi^{G,q}$ or $\Phi^{G,M}$.

\bibliographystyle{hplain}
\bibliography{database}

\end{document}